\documentclass[leqno,12pt]{article}
\usepackage{ams}
\usepackage{times}
\newtheorem{definition}{Definition}[section]
\newtheorem{lemma}[definition]{Lemma}
\newtheorem{prop}[definition]{Proposition}
\newtheorem{thm}[definition]{Theorem}
\newtheorem{cor}[definition]{Corollary}
\newtheorem{rmk}[definition]{Remark}

\def\newline{\hfil\break}

\def\dot{\hskip -.2cm {\bf .}\hskip .2cm}

\def\mapright#1{\mathop{\vbox{\ialign{
                ##\crcr
    ${\scriptstyle\hfil\;\;#1\;\;\hfil}$\crcr
 \noalign{\kern-1pt\nointerlineskip}
    \rightarrowfill\crcr}}\;}}

\def\mapleft#1{\mathop{\vbox{\ialign{
                ##\crcr
    ${\scriptstyle\hfil\;\;#1\;\;\hfil}$\crcr
    \noalign{\kern-1pt\nointerlineskip}
    \leftarrowfill\crcr}}\;}}

\def\proof{
  \noindent
  {\bf Proof:}
}
\def\proofth{
  \noindent
  {\bf Proof of Theorem }
}

\def\endproof{
{\unskip\nobreak\hfil\penalty50\hskip2em\hbox{}\nobreak\hfill
          $\square$\bigbreak}
}

\usepackage{amssymb}
\begin{document}

\title{Pseudo Harmonic Morphisms on Riemannian Polyhedra}
\author{M. A. Aprodu, T. Bouziane}
\date{}
\maketitle

\begin{abstract}
\noindent

The aim of this paper is to extend the notion of pseudo harmonic
morphism (introduced by Loubeau \cite {Lo}) to the case when the
source manifold is an admissible Riemannian polyhedron. We define
these maps to be harmonic in the sense of Eells-Fuglede \cite {EF}
and pseudo-horizontally weakly conformal in our sense (see Section
3). We characterize them by  means of  germs of harmonic functions
on the source polyhedron, in sense of Korevaar-Schoen \cite {KS},
and germs of holomorphic functions on the K\"ahler target
manifold.

\bigskip


\end{abstract}

\section{Introduction.}

"Harmonicity" is a topic which is situated between geometry and
analysis. For instance, Fuglede \cite{F} and Ishihara \cite{I},
independently, proved that harmonic morphisms between smooth
Riemannian manifolds (maps which pull back germs of harmonic
functions to germs of harmonic functions) are precisely harmonic
maps (analytic property) which are horizontally weakly conformal
(geometric property). A natural question arises: is there any
equivalent notion if the target manifolds are Hermitian or
K\"{a}hler?  If yes, can we characterize geometrically this
notion? Loubeau, in \cite{Lo}, gave  complete answers to these
questions and named the maps "Pseudo harmonic morphism".

In \cite{KS}, Korevaar and Schoen extended the theory of harmonic
maps between smooth Riemannian manifolds to the case of maps
between certain singular spaces: for example

\bigskip

\noindent\rule{12.1cm}{.2mm}

\noindent{\small M. A. Aprodu: Department of Mathematics,
University of Gala\c ti,
Domneasc\u a Str. 47, RO-6200, Gala\c ti, Romania.\\
e-mail: Monica.Aprodu@ugal.ro\\
T. Bouziane: The Abdus Salam International Center for Theoretical
Physics, strada Costiera 11, 34014 Trieste , Italy.\\
e-mail: tbouzian@ictp.trieste.it\\

\noindent{\em Mathematics Subject Classification (2000)}: 58E20,
53C43, 53C55, 32Q15. }
\newpage

\noindent admissible Riemannian polyhedra. The Riemannian
polyhedra are both very interesting examples of the "geometric
habitat" of the harmonicity (being harmonic spaces) and provide
several examples: smooth Riemannian manifolds, Riemannian orbit
spaces, normal analytic spaces, Thom spaces etc. Later, Eells and
Fuglede in \cite{EF}, expanded the notion of harmonic morphisms to
the case of Riemannian polyhedra. But, to give the same
characterization for harmonic morphisms between Riemannian
polyhedra as Fuglede and Ishihara did in the smooth case, they had
to pay a price:  the target had  to be a smooth Riemannian
manifold. Also, many of the properties found for the harmonic maps
and harmonic morphisms in the smooth case could be recovered when
consider as domain and target Riemannian polyhedra.

Remaining in the same spirt of ideas, the aim of this paper is to
extend
 pseudo harmonic morphisms to the case when the domain is an
admissible Riemannian polyhedron and the target a K\"ahler
manifold and to characterize them (as it was done in the smooth
case) by
 "geometric criteria" and "analytic criteria".
It turns out, because of the absence of global differential
calculus on singular spaces, that it is not easy to find a good
definition of the pseudo harmonic morphisms on Riemannian
polyhedra generalizing in a natural way the smooth case. Another
difficulty, compare with Loubeau's results, is to find a geometric
condition which characterize  pseudo harmonic morphisms on
Riemannian polyhedra, knowing that we can not talk about
horizontal vectors for example. A third difficulty is in the use
of  germs of harmonic functions in the sense of Korevaar-Schoen,
as the analytic aspect of our construction.

The outline of the paper is as follows.  In Section 2 we recall
some  results on Riemannian polyhedra, harmonic maps and morphisms
between Riemannian polyhedra. Section 3 is devoted to the
(alternative) geometric characterization of the pseudo harmonic
morphisms defined on Riemannian polyhedra named (also) "
pseudo-horizontally weak conformality ".  We  show that this
geometric property is preserved by the holomorphy. In Section 4,
we introduce the notion of pseudo harmonic morphisms from a
Riemannian polyhedra to a K\"ahler manifold and characterize them
as maps which pull back germs of holomorphic functions on target
manifold to germs of harmonic functions on the Riemannian
polyhedron (Theorem \ref{thm:01}). We also state a lifting
property for pseudo harmonic morphisms (Proposition
\ref{prop:fin}).

Finally, in Section 5, applying Proposition \ref{prop:fin}, we
give some interesting examples.

\section{Preliminaries.}

This section is devoted to some basic notions and known results
which will be used in the next sections.

\subsection{Riemannian polyhedra.}

\subsubsection{ Riemannian admissible complexes (\cite{BB}, \cite{Br}, \cite{BH},
 \cite{DJ}, \cite{T}).}

Let $K$ be a locally finite simplicial complex, endowed with a
 piecewise smooth Riemannian metric $g$ ( i.e. $g$ is a family of
 smooth Riemannian metrics $g_\Delta$ on simplices $\Delta$ of $K$,
 such that the restriction ${(g_\Delta)}_{|\Delta'}=g_{\Delta'}$, for any
 simplices $\Delta'$ and $\Delta$ with $\Delta'\subset \Delta$).

 Let $K$ be a finite dimensional simplicial complex which is connected
locally finite. A map $f$ from $[a,b]$ to $K$ is called a broken
geodesic if there is a subdivision $a=t_0<t_1<...<t_{p+1}=b$, such
that $f([t_i,t_{i+1}])$ is contained in some cell and the
restriction of $f$ to $[t_i,t_{i+1}]$ is a geodesic inside that
cell. Then define the length of the broken geodesic map $f$ to be:
 $$ L(f)=\sum_{i=0}^p d(f(t_i),f(t_{i+1})).$$
 The length inside each cell being measured with respect to its metric.

 Then, define $\tilde d(x,y)$, for every two points $x,y$ in $K$,
 to be the lower bound of the lengths of broken geodesics from $x$
 to $y$. $\tilde d$ is a pseudo-distance.

 If $K$ is connected and locally finite, then
$(K,\tilde d)$ is a length space and hence a geodesic space (i.e.
a metric space where every two points are connected by a curve
with length equal to the distance between them ) if complete.

A $l$-simplex in $K$ is called a {\it boundary simplex} if it is
adjacent to exactly one $l+1$ simplex. The complex $K$ is called
{\it boundaryless} if there are no boundary simplices in $K$.

The (open) {\it star} of an open simplex $\Delta^o$ (i.e. the
topological interior of $\Delta$ or the points of $\Delta$ not
belonging to any sub-face of $\Delta$; if $\Delta$ is point then
$\Delta^o=\Delta$) of $K$ is defined as:
$$
st(\Delta^o)=\bigcup \{ \Delta_i^o : \Delta_i \mbox{ is simplex of
} K \mbox { with } \Delta_i\supset\Delta \}.
$$
The star $st(p)$ of point $p$ is defined as the star of its {\it
carrier}, the unique open simplex $\Delta^o$ containing $p$. Every
star is path connected and contains the star of its points. In
particular $K$ is locally path connected. The closure of any star
is sub-complex.

\hskip 0 cm

 We say that the complex $K$ is {\it admissible}, if it is dimensionally homogeneous
and for every connected open subset $U$ of $K$, the open set
$U\setminus \{ U\cap \{\mbox{ the }(k-2)-\mbox{ skeleton }\} \}$
is connected, where $k$ is the dimension of $K$ (i.e. $K$ is
$(n-1)$-chainable).

\hskip 0 cm

Let $x\in K$ a vertex of $K$ so that $x$ is in the $l$-simplex
$\Delta_{l}$. We view $\Delta_{l}$ as an affine simplex in ${\bf
R}^l$, that is $\Delta _l =\bigcap_{i=0}^l H_i$, where
$H_0,H_1,...,H_l$ are closed half spaces in general position, and
we suppose that $x$ is in the topological interior of $H_0$. The
Riemannian metric $g_{\Delta_l}$ is the restriction to $\Delta_l$
of a smooth Riemannian metric defined in an open neighborhood $V$
of $\Delta_l$ in ${\bf R}^l$. The intersection
$T_x\Delta_l=\bigcap_{i=1}^l H_i \subset T_xV$ is a cone with apex
$0\in T_xV$, and $g_{\Delta_l}(x)$ turns it into an Euclidean
cone. Let $\Delta_m\subset \Delta_l$ ($m<l$) be another simplex
adjacent to $x$. Then, the face of $T_x\Delta_l$ corresponding to
$\Delta_m$ is isomorphic to $T_x\Delta_m$ and we view
$T_x\Delta_m$ as a subset of $T_x\Delta_l$.

Set $T_xK =\bigcup_{\Delta_i\ni x} T_x\Delta_i$ and we call it the
{\it tangent cone} of $K$ at $x$. Let $S_x\Delta_l$ denote the
subset of all unit vectors in $T_x\Delta_l$ and set $S_x=S_xK
=\bigcup_{\Delta_i\ni x} S_x\Delta_i$. The set $S_x$ is called the
{\it link} of $x$ in $K$. If $\Delta_l$ is a simplex adjacent to
$x$, then $g_{\Delta_l}(x)$ defines a Riemannian metric on the
$(l-1)$-simplex $S_x\Delta_l$. The family $g_x$ of riemannian
metrics $g_{\Delta_l}(x)$ turns $S_x\Delta_l$ into a simplicial
complex with a piecewise smooth Riemannian metric such that the
simplices are spherical.

We call an admissible  connected locally finite simplicial
complex, endowed with a piecewise smooth Riemannian metric, an
{\it admissible Riemannian complex}.

\hskip 0cm

\subsubsection{Riemannian polyhedron \cite{EF}.}

We mean by {\it polyhedron} a connected locally compact separable
Hausdorff space $X$ for which there exists a simplicial complex
$K$ and homeomorphism $\theta : K \rightarrow X$. Any such pair
$(K,\theta )$ is called a {\it triangulation} of $X$. The complex
$K$ is necessarily countable and locally finite (cf. \cite{S} page
120) and the space $X$ is path connected and locally contractible.
The {\it dimension} of $X$ is by definition the dimension of $K$
and it is independent of the triangulation.

A {\it sub-polyhedron} of a polyhedron $X$ with given
triangulation $(K,\theta )$, is the polyhedron $X'\subset X$
having as a triangulation $(K',\theta_{|K'})$ where $K'$ is a
subcomplex of $K$ (i.e. $K'$ is the complex whose vertices and
simplexes are some of those of $K$).

If $X$ is a polyhedron with specified triangulation $(K,\theta)$,
we shall speak of vertices, simplexes, $i$-skeletons or stars of
$X$ respectively of a space of links or tangent cones of $X$ as
the image under $\theta$ of vertices, simplexes, $i$-skeletons or
stars of $K$ respectively the image of space of links or tangent
cones of $K$. Thus our simplexes become compact subsets of $X$ and
the $i-$skeletons and stars become sub-polyhedrons of $X$.

If for given triangulation $(K,\theta)$ of the polyhedron $X$, the
homeomorphism $\theta$ is locally bi-lipschitz then $X$ is said
{\it Lip polyhedron} and $\theta$ {\it Lip homeomorphism}.

A {\it null set} in a Lip polyhedron $X$ is a set $Z\subset X$
such that $Z$ meets every maximal simplex $\Delta$, relative to a
triangulation $(K,\theta)$ (hence any) in set whose pre-image
under $\theta$ has $n$-dimensional Lebesgue measure $0$,
$n=dim\Delta$. Note that {\it 'almost everywhere' } (a.e.) means
everywhere except in some null set.

A {\it Riemannian polyhedron} $X=(X,g)$  is defined as a Lip
polyhedron $X$ with a specified triangulation $(K,\theta)$ such
that K is a simplicial complex endowed with a covariant bounded
measurable Riemannian metric tensor $g$, satisfying the
ellipticity condition below. In fact, suppose that $X$ has
homogeneous dimension $n$ and choose a measurable riemannian
metric $g_\Delta$ on the open euclidean $n$-simplex
$\theta^{-1}(\Delta^o)$ of $K$. In terms of euclidean coordinates
$\{x_1,...,x_n\}$ of points $x=\theta^{-1}(p)$, $g_\Delta$ thus
assigns to almost every point $p\in \Delta^o$ (or $x$), an
$n\times n$ symmetric positive definite matrix $g_\Delta =
(g_{ij}^\Delta(x))_{i,j=1,...,n}$ with measurable real entries and
there is a constant $\Lambda_\Delta >0$ such that (ellipticity
condition):
$$
\Lambda_\Delta^{-2}\sum_{i=0}^n(\xi^i)^2\leq \sum_{i,j}
g^\Delta_{ij}(x)\xi^i\xi^j\leq\Lambda_\Delta^2\sum_{i=0}^n(\xi^i)^2
$$
for $a.e.$ $x\in\theta^{-1}(\Delta^o)$ and every
$\xi=(\xi^1,...,\xi^n) \in {\bf R}^n$. This condition amounts to
the components of $g_\Delta$ being bounded and it is independent
not only of the choice of the euclidean frame on
$\theta^{-1}(\Delta^o)$ but also of the chosen triangulation.

For simplicity of statements we shall sometimes require that,
relative to a fixed triangulation $(K,\theta)$ of Riemannian
polyhedron $X$ (uniform ellipticity condition),
$$
\Lambda:= \mbox{ sup } \{\Lambda_\Delta:\Delta \mbox{ is simplex
of } X\}<\infty.
$$

A Riemannian polyhedron $X$ is said to be {\it admissible} if for
a fixed triangulation $(K,\theta)$ (hence any) the Riemannian
simplicial complex $K$ is admissible.

We underline that (for simplicity) the given definition of a
Riemannian polyhedron $(X,g)$ contains already the fact (because
of the definition above of the Riemannian admissible complex) that
the metric $g$ is {\it continuous} relative to some (hence any)
triangulation (i.e. for every maximal simplex $\Delta$ the metric
$g_\Delta$ is continuous up to the boundary). This fact is
sometimes in the literature omitted. The polyhedron is said to be
{\it simplexwise smooth } if relative to some triangulation
$(K,\theta)$ (and hence any), the complex $K$ is simplexwise
smooth. Both continuity and simplexwise smoothness are preserved
under subdivision.

In the case of a general bounded measurable Riemannian metric $g$
on $X$, we often consider, in addition to $g$, the {\it euclidean
Riemannian metric} $g^e$ on the Lip polyhedron $X$ with a
specified triangulation $(K,\theta)$. For each simplex $\Delta$,
$g^e_\Delta$ is defined in terms of euclidean frame on
$\theta^{-1}(\Delta^o)$ as above by unit matrix
$(\delta_{ij})_{i,j}$. Thus $g^e$ is by no means covariantly
defined and should be regarded as a mere reference metric on the
triangulated polyhedron $X$.

Relative to a given triangulation $(K,\theta)$ of an
$n$-dimensional Riemannian polyhedron $(X,g)$ (not necessarily
admissible), we have on $X$ the distance function $e$ induced by
the euclidean distance on the euclidean space $V$ in which $K$ is
affinely Lip embedded. This distance $e$ is not intrinsic but it
will play an auxiliary role in defining an equivalent distance
$d_X$ as follows:

Let $\frak Z$ denote the collection of all null sets of $X$. For
given triangulation $(K,\theta)$ consider the set $Z_K\subset
\frak Z $ obtained from $X$ by removing from each maximal simplex
$\Delta$ in $X$ those points of $\Delta^o$ which are Lebesgue
points for $g_\Delta$. For $x,y \in X$ and any $Z\in \frak Z$ such
that $Z\subset Z_K$ we set:
$$
d_X(x,y)
  =\sup\limits_{\begin{array}{cc}
                            Z\in \frak Z \\
                            Z\supset Z_K\\
                      \end{array}}
         \inf\limits_{\begin{array}{ccc}
                             \gamma\\
                             \gamma(a)=x,\\
                              \gamma(b)=y\\
                      \end{array}}
 \left\{ L_K(\gamma):
                       \begin{array}{ll}
                        \gamma & \mbox{is Lip continuous path}\\
                                &\mbox{ and transversal to } Z\\
                       \end{array}
 \right\},
$$
where $L_K(\gamma)$ is  the length of the path $\gamma$ defined
as:
$$
L_K(\gamma)= \sum\limits_{\Delta\subset X}
 \int\limits_{\gamma^{-1}(\Delta^o)}
 \sqrt{(g_{ij}^\Delta \circ \theta^{-1} \circ \gamma)
 \dot{\gamma}^i \dot{\gamma}^j } ,
 \begin{array}{l}
 \mbox{ the sum is over}\\
 \mbox{ all simplexes meeting } \gamma.\\
 \end{array}
 $$

 It is shown in \cite{EF} that the distance $d_X$ is
 intrinsic, in particular it is independent of the chosen triangulation
 and it is equivalent to the euclidean distance $e$ (due to the Lip
 affinely and homeomorphically embedding of $X$ in some euclidean space $V$).

\subsection{Energy of maps}

The concept of energy in the case of a map of Riemannian domain
into an arbitrary metric space $Y$ was defined and investigated by
Korevaar and Schoen \cite{KS}. Later this concept was extended by
Eells and Fuglede \cite{EF} to the case of map from an admissible
Riemannian polyhedron $X$ with simplexwise smooth Riemannian
metric. Thus, the {\it energy } $E(\varphi)$ of a map $\varphi$
from $X$ to the space $Y$ is defined as the limit of suitable
approximate energy expressed in terms of the distance function
$d_Y$ of $Y$.

 It is shown in \cite{EF} that the maps $\varphi : X\rightarrow Y$ of
 finite energy are precisely those quasicontinuous (i.e.
 has a continuous restriction to closed sets), whose complements
 have arbitrarily small capacity, (cf. \cite{EF} page 153) whose restriction
 to each top dimensional simplex of $X$
 has finite energy in the sense of Korevaar-Schoen, and
 $E(\varphi)$ is the sum of the energies of these restrictions.

 Now, let $(X,g)$ be an admissible $m$-dimensional
 Riemannian polyhedron with simplexwise smooth Riemannian metric.
 It is not required that $g$ is continuous across lower
 dimensional simplexes. The target $(Y,d_Y)$ is an arbitrary
 metric space.

Denote $L^2_{loc}(X,Y)$ the space of all $\mu_g$-measurable
 ($\mu_g$ the volume measure of $g$)
 maps $\varphi :X\rightarrow Y$ having separable essential range and
 for which the map $d_Y(\varphi (.),q)\in L^2_{loc}(X,\mu_g)$
 (i.e. locally $\mu_g$-squared integrable) for some point $q$
 (hence by triangle inequality for any point). For $\varphi,\psi \in
 L^2_{loc}(X,Y)$ define their distance $D(\varphi,\psi)$ by:
 $$
 D^2(\varphi,\psi)= \int\limits_X d_Y^2(\varphi (x),\psi(y)) d\mu_g(x).
 $$
 Two maps $\varphi,\psi \in L^2_{loc}(X,Y)$ are said to be {\it
 equivalent} if $D(\varphi,\psi)=0$,(i.e. $\varphi(x)=\psi(x)$
 $\mu_g$-a.e.). If the space $X$ is compact then
 $D(\varphi,\psi)<\infty$
 and $D$ is a metric on $L^2_{loc}(X,Y)=L^2(X,Y)$ and complete if
 the space $Y$ is complete \cite{KS}.

 The {\it approximate energy density} of the map $\varphi\in
 L^2_{loc}(X,Y)$ is defined for $\epsilon >0$ by:
 $$
 e_\epsilon(\varphi)(x)=\int\limits_{B_X(x,\epsilon)}\frac{d_Y^2(\varphi(x),\varphi(x'))}
 {\epsilon^{m+2}}d\mu_g(x').
 $$
 The function $e_\epsilon(\varphi)\ge 0$ is locally
 $\mu_g$-integrable.

The {\it energy} $E(\varphi)$ of a map $\varphi$ of class
 $L^2_{loc}(X,Y)$ is:
 $$
 E(\varphi)=\sup_{f\in
 C_c(X,[0,1])}(\limsup_{\epsilon\rightarrow 0}\int\limits_X
 fe_\epsilon(\varphi) d\mu_g),
 $$
 where $C_c(X,[0,1])$ denotes the
 space of continuous functions from $X$ to the interval $[0,1]$
 with compact support.

A map $\varphi: X\rightarrow Y$ is said to be {\it locally of
finite energy}, and we write $\varphi \in W^{1,2}_{loc}(X,Y)$, if
$E(\varphi_{|U})<\infty$ for every relatively compact domain
$U\subset X$, or equivalently if $X$ can be covered by domains
$U\subset X$ such that $E(\varphi_{|U})<\infty$.

For example (cf. \cite{EF} Lemma 4.4), every Lip continuous map
$\varphi: X \rightarrow Y$ is of class $W^{1,2}_{loc}(X,Y)$. In
the case when $X$ is compact $W^{1,2}_{loc}(X,Y)$ is denoted
$W^{1,2}(X,Y)$ the space of all maps of finite energy.

$W^{1,2}_c(X,Y)$ denotes the linear subspace of $W^{1,2}(X,Y)$
consisting of all maps of finite energy of compact support in $X$.

 We can show  (cf. \cite{EF} Theorem 9.1) that a real function
$\varphi \in L^2_{loc}(X)$ is locally of finite energy if and only
if  there is a function $e(\varphi)\in L^1_{loc}(X)$, named {\it
energy density} of $\varphi$, such that (weak convergence):
$$
\lim_{\epsilon\rightarrow 0}\int\limits_X f e_\epsilon (\varphi)
d\mu_g =\int\limits_X f e(\varphi) d\mu_g, \mbox{  for each  }
f\in C_c(X).
$$

\subsection{Harmonic maps and harmonic morphisms on Riemannian polyhedra \cite{EF}.}

In this paragraph we shall remind some relevant results which give
the relation between harmonic morphisms and harmonic maps on
Riemannian polyhedra.

\subsubsection{Harmonic maps.}

Let $(X,g)$ be an arbitrary admissible Riemannian polyhedron ($g$
just bounded measurable with local elliptic bounds), $dim X=m$ and
$(Y,d_Y)$ a metric space .

 A continuous map $\varphi:X\rightarrow Y$ of class $W_{loc}^{1,2}(X,Y)$ is said to be
 {\it harmonic} if it is {\it bi-locally E-minimizing}, i.e. $X$ can be covered by
relatively compact subdomains $U$ for each of which there is an
open set $V\supset \varphi(U)$ in $Y$ such that
$$
E(\varphi_{|U})\leq E(\psi_{|U})
$$
for every continuous map $\psi\in W_{loc}^{1,2}(X,Y)$, with
$\psi(U)\subset V$ and $\psi=\varphi$ in $X\backslash U$.

Let $(N,h)$ denote a smooth Riemannian manifold without boundary,
$dim_{\bf R} N=n$ and $\Gamma^k_{\alpha\beta}$ the Christoffel
symbols on $N$. By a {\it weakly harmonic map }
$\varphi:X\rightarrow N$ we mean  a quasicontinuous map (a map
which is continuous on the complement of open sets of arbitrarily
small capacity ; in the case of the Riemannian polyhedron $X$ it
is just the complement of open subsets of the $(m-2)$-skeleton of
$X$) of class $W_{loc}^{1,2}(X,N)$ with the following property:

For any chart $\eta:V\rightarrow {\bf R}^n$ on N and any quasiopen
set $U\subset\varphi^{-1}(V)$ of compact closure in X, the
equation
$$
\int\limits_U\langle\nabla\lambda,\nabla\varphi^k\rangle d\mu_g=
\int\limits_U \lambda(\Gamma^k_{\alpha\beta}\circ\varphi)
\langle\nabla\varphi^\alpha,\nabla\varphi^\beta\rangle d\mu_g,
$$
holds for every $k=1,...,n$ and every bounded function $\lambda\in
W_0^{1,2}(U)$.

It is shown in \cite{EF}, (Theorem 12.1), that: for a continuous
map $\varphi\in W_{loc}^{1,2}(X,N)$ the following are equivalent:

(a) $\varphi$ is harmonic,

(b) $\varphi$ is weakly harmonic,

(c) $\varphi$ pulls convex functions on open sets $V\subset N$
back to subharmonic functions on $\varphi^{-1}(V)$.

\subsubsection{Harmonic morphisms.}

Denote by $X$ and $Y$ two Riemannian polyhedra (or any harmonic
spaces in the sense of Brelot; see Chapter 2, \cite{EF}).

 A continuous map $\varphi:X\rightarrow Y$ is a {\it harmonic
morphism} if, for every open set $V\subset Y$ and for every
harmonic function $v$ on $V$, $v\circ\varphi$ is harmonic on
$\varphi^{-1}(V)$.

 Let $\varphi:X\rightarrow Y$ be a nonconstant
harmonic morphism, then (cf. \cite{EF}, Theorem 13.1):

(i) $\varphi$ likewise pulls germs of superharmonic functions on
$Y$ back to germs of superharmonic functions on $X$.

(ii) If $\varphi$ is surjective and proper then a function
$v:V\rightarrow [-\infty,\infty]$ ($V$ open in $Y$) is
superharmonic [resp. harmonic] if (and only if) $v\circ\varphi$ is
superharmonic [resp. harmonic] on $\varphi^{-1}(V)$.

 Let $(N,g_N)$ denote a $n$-Riemannian manifold without
boundary and suppose that the polyhedron $X$ is admissible. A
continuous map $\varphi:X\rightarrow N$ of class
$W^{1,2}_{loc}(X,N)$ is called {\it horizontally weakly conformal}
if there exist a scalar $\lambda$, defined $a.e.$ in $X$, such
that:
$$
\langle\nabla(v\circ\varphi),\nabla(w\circ\varphi)\rangle=
\lambda[g_N(\nabla_Nv,\nabla_Nw)\circ\varphi] \mbox{ a.e. in X }
$$
for every pair of functions $v,w\in \mathcal{C}^1(N)$. Henceforth
$\nabla_N$ denote the gradient operator on $N$ and $\nabla$ the
gradient operator defined a.e. on the domain space $(X,g)$.

The property of horizontally weak conformality is a local one,
thus it reads in terms of local coordinates $(y^\alpha)$ in $N$,

$$
\langle\nabla\varphi^\alpha,\nabla\varphi^\beta\rangle=
\lambda(g_N^{\alpha\beta}\circ\varphi) \mbox{ a.e. in X }
$$
for $\alpha,\beta=1,...,n$. Taking $\alpha=\beta$, $\lambda$ is
uniquely determined and $\lambda\geq 0$ a.e. in $X$. Moreover,
$\lambda\in L^1_{loc}(X)$ because $\nabla\varphi^\alpha\in
L^1_{loc}(X)$. $\lambda$ is called the {\em dilation} of
$\varphi$.

The notion of horizontally weak conformality is intimately related
to the one of harmonic morphisms. For instance we can show (cf.
\cite{EF}, Theorem 13.2) that  a continuous map
$\varphi:X\rightarrow N$ of class $W^{1,2}_{loc}(X,N)$ is a
harmonic morphism if and only if $\varphi$ is horizontally weakly
conformal, harmonic map and equivalently, there is a scalar
$\lambda\in L^1_{loc}(X)$ such that
$$
-\int\limits_X\langle\nabla\psi,\nabla(v\circ\varphi)\rangle=
\int\limits_X\psi\lambda[(\Delta_Nv)\circ\varphi]
$$
for every $v\in \mathcal{C}^2(N)$ and $\psi\in Lip_c(X)$ ( or
$\psi\in W^{1,2}_0(X)\cap L^\infty(X)$).

In the affirmative case, $\lambda$ from the last equality equals
a.e. the dilation of $\varphi$ (as a horizontally weakly conformal
map).

As a consequence, for a harmonic morphism $\varphi:X\rightarrow
N$, if $\psi:N\rightarrow P$ is a harmonic map between smooth
Riemannian manifolds without boundary, then the composition
$\psi\circ\varphi$ is a harmonic map.

\section{Pseudo-horizontally weakly conformal maps on Riemannian polyhedra.}

The aim of the present section, is to extend the notion of {\it
pseudo-horizontally weakly conformal} maps on Riemannian manifolds
(see \cite{Lo}) to Riemannian polyhedra and to establish their
properties. We will use the same terminology as in \cite{Lo}.

Let $(X,g)$ be an admissible Riemannian polyhedron of $dim X=m$
and $(N,J^N,g_N)$ a Hermitian manifold of ${dim_{\bf R}}N=2n$,
without boundary.

We  denote by  $Holom (N)=\{f:N\rightarrow {\bf C}, f  \mbox {
local holomorphic function}\}$. In what follows, the gradient
operator  and the inner product in $(X,g)$ are well defined a.e.
in $X$ and will be denoted by $\nabla$ and $\langle , \rangle$
respectively.

\begin{definition}\label{def:01}
Let $\varphi:X\rightarrow N$ be a continuous map of class
$W_{loc}^{1,2}(X,N)$. $\varphi$ is called {\rm pseudo-horizontally
weakly conformal} (shortening PHWC), if for any pair of local
holomorphic functions $v,w\in Holom (N)$, such that $v=v_1+iv_2$,
$w=w_1+iw_2$, we have:
\begin{equation}\label{01}
\left\{
\begin{array}{cc}
\langle\nabla(w_1\circ\varphi), \nabla(v_1\circ\varphi)\rangle-
\langle\nabla(w_2\circ\varphi), \nabla(v_2\circ\varphi)\rangle=0& \mbox{  a.e. in }X \\
\langle\nabla(w_2\circ\varphi), \nabla(v_1\circ\varphi)\rangle+
\langle\nabla(w_1\circ\varphi), \nabla(v_2\circ\varphi)\rangle=0& \mbox{  a.e. in }X \\
\end{array}
\right.
\end{equation}
\end{definition}

\begin{rmk}\label{rmk:01}
 {\rm Definition \ref {def:01}} is a local one, hence it is
sufficient to check the identities (\ref {01}) in local complex
coordinates $(z_1,z_2,...,z_n)$ in $N$. Taking $z_A=x_A+iy_A,
\forall A=1,...,n$, the relations (\ref {01}), $\forall
A,B=1,...,n$, read:
\begin{equation}\label{02}
\left\{
\begin{array}{cc}
\langle\nabla\varphi_1^B, \nabla\varphi_1^A\rangle-
\langle\nabla\varphi_2^B, \nabla\varphi_2^A\rangle=0& \mbox{  a.e. in }X \\
&\\
\langle\nabla\varphi_2^B, \nabla\varphi_1^A\rangle+
\langle\nabla\varphi_1^B, \nabla\varphi_2^A\rangle=0& \mbox{  a.e. in }X \\
\end{array}
\right.
\end{equation}
where $$ \left\{
\begin{array}{cc}
\varphi_1^A:=x_A\circ\varphi,\mbox{  }
\varphi_2^A:=y_A\circ\varphi, &\forall A=1,...,n ,\\
&\\
\varphi_1^B:=x_B\circ\varphi,\mbox{  }
\varphi_2^B:=y_B\circ\varphi, &\forall B=1,...,n.\\
\end{array}
\right.
$$
\end{rmk}

\begin{rmk}\label{rmk:def}
{\rm Definition \ref{def:01}} is justified by seeing that if the
source manifold is a smooth Riemannian one, without boundary, we
obtain exactly the commuting condition between $d\varphi_x\circ
d\varphi^*_x$ and $J^N_{\varphi(x)}$ {\rm (see \cite {Lo} or \cite
{BW}, \cite {AAB})}, where $d\varphi^*_x:T_{\varphi(x)}N
\rightarrow T_{x}X$ is the adjoint map of the tangent linear map
$d\varphi_x:T_{x}X\rightarrow T_{\varphi(x)}N$, for any $x\in X$.
\end{rmk}

The next proposition justifies the use of the term 'horizontally
weakly conformal', indeed we obtain, when the target dimension is
one, an equivalence between the horizontally weakly conformality m
and pseudo-horizontally weakly conformality, as in the smooth
case.

\begin{prop}\label{prop:00}
Let $\varphi:X\rightarrow N$ a horizontally weakly conformal map
{\rm (see subsection 2.3.2)} from a Riemannian admissible
polyhedron $(X,g)$ into a Hermitian manifold $(N,J^N,g_N)$. Then
$\varphi$ is pseudo-horizontally weakly conformal. If the complex
dimension of $N$ is equal to one, then the two conditions are
equivalent.
\end{prop}

\proof Let $\varphi:X\rightarrow N$ be a horizontally weakly
conformal map from an admissible Riemannian polyhedron $(X,g)$
into a Hermitian manifold $(N,J^N,g_N)$ of real dimension $2n$.
Take $(z_A=x_A+iy_A)_{A=1,...,n}$ local complex coordinates in
$N$. Then $\big\{\frac{\partial}{\partial x_1},...,
\frac{\partial}{\partial x_n},\frac{\partial}{\partial y_1},...,
\frac{\partial}{\partial y_n}\big\}$
 is a local frame in $TN$ such that
 $$
 \left\{
 \begin{array}{c}
   J^N(\frac{\partial}{\partial x_A})=
\frac{\partial}{\partial y_A} \\
\\
J^N(\frac{\partial}{\partial y_A})=-
\frac{\partial}{\partial x_A} \\
 \end{array}
 \right.,
 \forall A=1,...,n.
 $$
The horizontally weakly conformal condition reads in the
considered frame:
\begin{equation}\label{04}
\langle\nabla\varphi^\alpha,\nabla\varphi^\beta\rangle=
\lambda(g_N^{\alpha\beta}\circ\varphi) \mbox{  a.e. in  }X,
\forall\alpha,\beta=1,...,2n.
\end{equation}
where $\varphi^\alpha=\xi_\alpha\circ\varphi$ for $ \xi_\alpha=
\left\{
\begin{array}{cc}
  x_\alpha ,& \alpha=1,...,n \\
  y_{\alpha-n} ,& \alpha=n+1,...,2n\\
\end{array}
\right.. $ Explicitly the equalities (\ref {04}) are the
following:

\begin{equation}\label{05}
\left\{
\begin{array}{cc}
 \langle\nabla(\xi_\alpha\circ\varphi),\nabla(\xi_\beta\circ\varphi)\rangle=
\lambda[g_N(\frac{\partial}{\partial\xi_\alpha},
\frac{\partial}{\partial\xi_\beta})\circ\varphi]& \mbox{ a.e. in  }X,\\
      &\forall\alpha,\beta=1,...,n.\\
  \langle\nabla(\xi_\alpha\circ\varphi),\nabla(\xi_\beta\circ\varphi)\rangle=
\lambda[g_N(\frac{\partial}{\partial\xi_\alpha},
\frac{\partial}{\partial\xi_\beta})\circ\varphi] &\mbox{ a.e. in  }X,\\
     &\forall\alpha,\beta=n+1,...,2n. \\
\end{array}
\right.
\end{equation}
and
\begin{equation}\label{06}
\left\{
\begin{array}{cc}
 \langle\nabla(\xi_\alpha\circ\varphi),\nabla(\xi_\beta\circ\varphi)\rangle=
\lambda[g_N(\frac{\partial}{\partial\xi_\alpha},
\frac{\partial}{\partial\xi_\beta})\circ\varphi]& \mbox{ a.e. in  }X,\\
      &\forall\alpha=1,...,n;\\
      & \forall\beta=n+1,...,2n.\\
  \langle\nabla(\xi_\alpha\circ\varphi),\nabla(\xi_\beta\circ\varphi)\rangle=
\lambda[g_N(\frac{\partial}{\partial\xi_\alpha},
\frac{\partial}{\partial\xi_\beta})\circ\varphi] &\mbox{ a.e. in  }X,\\
     &\forall\alpha=n+1,...,2n;\\
     & \forall\beta=1,...,n. \\
\end{array}
\right.
\end{equation}
Read (\ref{05}) and (\ref{06}) in terms of $x_A$ and $y_B$,
$\forall A,B=1,...,n$:
\begin{equation}\label{07}
\left\{
\begin{array}{cc}
 \langle\nabla(x_A\circ\varphi),\nabla(x_B\circ\varphi)\rangle=
\lambda[g_N(\frac{\partial}{\partial x_A},
\frac{\partial}{\partial x_B})\circ\varphi]& \mbox{ a.e. in  }X,\\
&\\
\langle\nabla(y_A\circ\varphi),\nabla(y_B\circ\varphi)\rangle=
\lambda[g_N(\frac{\partial}{\partial y_A},
\frac{\partial}{\partial y_B})\circ\varphi] &\mbox{ a.e. in  }X.\\
\end{array}
\right.
\end{equation}
and
\begin{equation}\label{08}
\left\{
\begin{array}{cc}
 \langle\nabla(x_A\circ\varphi),\nabla(y_B\circ\varphi)\rangle=
\lambda[g_N(\frac{\partial}{\partial x_A},
\frac{\partial}{\partial y_B})\circ\varphi]& \mbox{ a.e. in  }X,\\
&\\
\langle\nabla(y_A\circ\varphi),\nabla(x_B\circ\varphi)\rangle=
\lambda[g_N(\frac{\partial}{\partial y_A},
\frac{\partial}{\partial x_B})\circ\varphi] &\mbox{ a.e. in  }X.\\
\end{array}
\right.
\end{equation}
Because $J^N$ is the complex structure with respect to the
hermitian metric $g_N$, we have $\forall A,B=1,...,n$:
\begin{equation}\label{09}
\left\{
\begin{array}{c}
g_N(\frac{\partial}{\partial x_A}, \frac{\partial}{\partial
x_B})=g_N(\frac{\partial}{\partial y_A},
\frac{\partial}{\partial y_B}),\\
\\
g_N(\frac{\partial}{\partial x_A}, \frac{\partial}{\partial
y_B})=-g_N(\frac{\partial}{\partial y_A},
\frac{\partial}{\partial x_B}).\\
\end{array}
\right.
\end{equation}
Invoking (\ref{09}), (\ref{07}) and (\ref{08}) we conclude that
$\varphi$ is pseudo-horizontally weakly conformal.

Consider now the case when $dim_{\bf C}N=1$ and suppose
$\varphi:X\rightarrow N$ is a pseudo-horizontally weakly conformal
map.

Let $z=x+iy$ be  a local complex chart in $N$. In terms of this
chart the pseudo-horizontally weakly conformal condition (\ref
{01}) reads:
\begin{equation}\label{10}
\left\{
\begin{array}{cc}
\langle\nabla\varphi_x,\nabla\varphi_x\rangle-
\langle\nabla\varphi_y,\nabla\varphi_y\rangle=0&
\mbox {a.e. in } X\\
\langle\nabla\varphi_y,\nabla\varphi_x\rangle+
\langle\nabla\varphi_x,\nabla\varphi_y\rangle=0&
\mbox {a.e. in } X\\
\end{array}
\right.
\end{equation}
Remember that $ g_N(\frac{\partial}{\partial x},
\frac{\partial}{\partial x})=g_N(\frac{\partial}{\partial y},
\frac{\partial}{\partial y})\neq 0 $, so we can define
\begin{equation}\label{11}
\lambda:=\frac{\langle\nabla\varphi_x,\nabla\varphi_x\rangle}
{g_N(\frac{\partial}{\partial x}, \frac{\partial}{\partial
x})\circ \varphi}=
\frac{\langle\nabla\varphi_y,\nabla\varphi_y\rangle}
{g_N(\frac{\partial}{\partial y}, \frac{\partial}{\partial
y})\circ \varphi}, \mbox{ a.e. in } X.
\end{equation}
From (\ref{10}) and (\ref{11}) we get:
\begin{equation}\label{12}
\left\{
\begin{array}{cc}
\langle\nabla\varphi_x,\nabla\varphi_x\rangle=\lambda[g_N(\frac{\partial}{\partial
x},
\frac{\partial}{\partial x})\circ \varphi]& \mbox{ a.e. in } X\\
&\\
\langle\nabla\varphi_y,\nabla\varphi_y\rangle=\lambda
[g_N(\frac{\partial}{\partial y},
\frac{\partial}{\partial y})\circ \varphi]& \mbox{ a.e. in } X\\
&\\
\langle\nabla\varphi_x,\nabla\varphi_y\rangle=0& \mbox {a.e. in } X\\
\end{array}
\right.
\end{equation}
which means that $\varphi$ is horizontally weakly conformal.
\endproof

PHWC maps on Riemannian polyhedra can be characterized thanks to
germs of  holomorphic functions on the target Hermitian manifolds
as follows:

\begin{prop}\label{prop:01'}
Let $\varphi:X\rightarrow N$ be a continuous map of class
$W_{loc}^{1,2}(X,N)$. Then $\varphi$ is pseudo-horizontally weakly
conformal if and only if for any local holomorphic function
$\psi:N\rightarrow {\bf C}$, $\psi\circ\varphi$ is also
pseudo-horizontally weakly conformal.
\end{prop}

\proof Let $\varphi:X\rightarrow N$ be a continuous map of class
$W_{loc}^{1,2}(X,N)$, and $\psi:N\rightarrow {\bf C}$ be any
holomorphic function with $\psi=\psi_1+ i \psi_2$.

It is obvious (by definition) that if $\varphi$ is a PHWC map then
$\psi\circ \varphi$ is a PHWC map.

Conversely, suppose now that for any holomorphic function
$\psi:N\rightarrow {\bf C}$, the composition
$\psi\circ\varphi:X\rightarrow {\bf C}$ is a PHWC function.
Throughout the proof of Proposition \ref{prop:00} we have seen
that this last fact reads:

\begin{equation}\label{12}
\left\{
\begin{array}{cc}
\langle\nabla(\psi_1\circ\varphi),\nabla(\psi_1\circ\varphi)\rangle-
\langle\nabla(\psi_2\circ\varphi),\nabla(\psi_2\circ\varphi)\rangle=0&
\mbox {  a.e. in } X\\
\langle\nabla(\psi_1\circ\varphi),\nabla(\psi_2\circ\varphi)\rangle+
\langle\nabla(\psi_2\circ\varphi),\nabla(\psi_1\circ\varphi)\rangle
=0&
\mbox {  a.e. in } X\\
\end{array}
\right.
\end{equation}

Then for a fixed local holomorphic chart
$(z_\alpha)_{\alpha=1,...,n}$ of $N$,
$z_\alpha=x_\alpha+iy_\alpha$, the equalities (\ref{12}) become:

\begin{equation}\label{13}
\left\{
\begin{array}{cc}
0=
\langle\nabla(\psi_1\circ\varphi),\nabla(\psi_1\circ\varphi)\rangle-
\langle\nabla(\psi_2\circ\varphi),\nabla(\psi_2\circ\varphi)\rangle=\\
\sum\limits_{k,l=1}^n\langle(\partial_k\psi_1\circ\varphi).\nabla\varphi^k,
(\partial_l\psi_1\circ\varphi).\nabla\varphi^l\rangle-\\
\sum\limits_{k,l=1}^n\langle(\partial_k\psi_2\circ\varphi).\nabla\varphi^k,
(\partial_l\psi_2\circ\varphi).\nabla\varphi^l\rangle &
\mbox {  a.e. in } X\\
&\\
\textrm{and}\\
&\\
0=
\langle\nabla(\psi_1\circ\varphi),\nabla(\psi_2\circ\varphi)\rangle+
\langle\nabla(\psi_2\circ\varphi),\nabla(\psi_1\circ\varphi)\rangle=\\
\sum\limits_{k,l=1}^n\langle(\partial_k\psi_1\circ\varphi).\nabla\varphi^k,
(\partial_l\psi_2\circ\varphi).\nabla\varphi^l\rangle +\\
\sum\limits_{k,l=1}^n\langle(\partial_k\psi_2\circ\varphi).\nabla\varphi^k,
(\partial_l\psi_1\circ\varphi).\nabla\varphi^l\rangle &
\mbox {  a.e. in } X\\
\end{array}
\right.
\end{equation}

\bigskip

where $\varphi^\alpha=\xi_\alpha\circ\varphi$ for $ \xi_\alpha=
\left\{
\begin{array}{cc}
  x_\alpha ,& \alpha=1,...,n \\
  y_{\alpha-n} ,& \alpha=n+1,...,2n\\
\end{array}
\right. .$

\bigskip

Taking into account the Cauchy-Riemann equations and after doing
some computations, (13) becomes:

\begin{equation}
\label{14} \left\{
\begin{array}{cc}
\sum\limits_{k,l=1}^n \left[ (\frac{\partial\psi_1}{\partial
x_k}\circ\varphi) (\frac{\partial\psi_1}{\partial
x_l}\circ\varphi)- (\frac{\partial\psi_2}{\partial
x_k}\circ\varphi) (\frac{\partial\psi_2}{\partial
x_l}\circ\varphi)
\right]& \\
\left[
\langle\nabla(x_k\circ\varphi),\nabla(x_l\circ\varphi)\rangle-
\langle\nabla(y_k\circ\varphi),\nabla(y_l\circ\varphi)\rangle
\right] &\\
+ &\\
\sum\limits_{k,l=1}^n \left[ (\frac{\partial \psi_1}{\partial
x_k}\circ\varphi) (\frac{\partial \psi_1}{\partial
y_l}\circ\varphi) - (\frac{\partial \psi_2}{\partial
x_k}\circ\varphi) (\frac{\partial \psi_2}{\partial
y_l}\circ\varphi)
\right]&\\
\left[ \langle\nabla(x_k\circ\varphi),
\nabla(y_l\circ\varphi)\rangle+ \langle\nabla(y_k\circ\varphi),
\nabla(x_l\circ\varphi)\rangle \right] =
0 &\mbox {   a.e. in } X \\
&\\
 \textrm{and} \\
 &\\
\sum\limits_{k,l=1}^n \left[ (\frac{\partial \psi_1}{\partial
x_k}\circ\varphi) (\frac{\partial \psi_2}{\partial
x_l}\circ\varphi) + (\frac{\partial \psi_2}{\partial
x_k}\circ\varphi) (\frac{\partial \psi_1}{\partial
x_l}\circ\varphi)
\right]&\\
\left[ \langle\nabla(x_k\circ\varphi),
\nabla(x_l\circ\varphi)\rangle- \langle\nabla(y_k\circ\varphi),
\nabla(y_l\circ\varphi)\rangle
\right]&\\
+&\\
\sum\limits_{k,l=1}^n \left[ (\frac{\partial \psi_1}{\partial
x_k}\circ\varphi) (\frac{\partial \psi_2}{\partial
y_l}\circ\varphi) + (\frac{\partial \psi_2}{\partial
x_k}\circ\varphi) (\frac{\partial \psi_1}{\partial
y_l}\circ\varphi)
\right]&\\
\left[ \langle\nabla(x_k\circ\varphi),
\nabla(y_l\circ\varphi)\rangle+ \langle\nabla(y_k\circ\varphi),
\nabla(x_l\circ\varphi)\rangle \right] =
0 &\mbox {   a.e. in } X\\
\end{array}
\right.
\end{equation}

Now, choose particular holomorphic functions $\psi$'s, for example
locally $\psi=z_k+z_l$ and vary the $k,l$ indices, we obtain all
the pseudo-horizontally weakly conformality conditions of the map
$\varphi$. This ends the proof.

\endproof

The next proposition makes clear the relation between  PHWC maps
on Riemannian polyhedra and  holomorphic maps on target Hermitian
manifolds.

\begin{prop}\label{prop:01}
Let $\varphi:X\rightarrow N$ a continuous map of class
$W_{loc}^{1,2}(X,N)$ and $(P,J^P,g_P)$ another Hermitian manifold
of $\rm{dim}_{\bf R}P=2p$. Then $\varphi$ is pseudo-horizontally
weakly conformal if and only if for every local holomorphic map
$\psi:N\rightarrow P$, $\psi\circ\varphi$ is also
pseudo-horizontally weakly conformal.
\end{prop}

\proof Let $\psi:N\rightarrow P$ be a local holomorphic map.
Choose $(z_\alpha)_{\alpha=1,...,p}$ local complex coordinates in
$P$ and denote $\psi^\alpha:=z_\alpha\circ\psi, \forall
\alpha=1,...,p$.

Suppose that $\varphi$ is pseudo-horizontally weakly conformal.
Then, by definition we have, for every pair of local holomorphic
functions $v,w\in Holom (N)$, such that $v=v_1+iv_2$,
$w=w_1+iw_2$,
\begin{equation}\label{}
\left\{
\begin{array}{cc}
\langle\nabla(w_1\circ\varphi), \nabla(v_1\circ\varphi)\rangle-
\langle\nabla(w_2\circ\varphi), \nabla(v_2\circ\varphi)\rangle=0& \mbox{ a.e. in }X \\
\langle\nabla(w_2\circ\varphi), \nabla(v_1\circ\varphi)\rangle+
\langle\nabla(w_1\circ\varphi), \nabla(v_2\circ\varphi)\rangle=0& \mbox{ a.e. in }X \\
\end{array}
\right.
\end{equation}
In particular, for every pair of local holomorphic functions on
the Hermitian manifold $N$, $\psi^\alpha$, $\psi^\beta$,
$\alpha,\beta = 1,...,p$ we have:

\begin{equation}\label{}
\left\{
\begin{array}{cc}
\langle\nabla(\psi^\alpha_1\circ\varphi),
\nabla(\psi^\beta_1\circ\varphi)\rangle-
\langle\nabla(\psi^\alpha_2\circ\varphi), \nabla(\psi^\beta_2\circ\varphi)\rangle=0& \mbox{ a.e. in }X \\
\langle\nabla(\psi^\alpha_2\circ\varphi),
\nabla(\psi^\beta_1\circ\varphi)\rangle+
\langle\nabla(\psi^\alpha_1\circ\varphi), \nabla(\psi^\beta_2\circ\varphi)\rangle=0& \mbox{ a.e. in }X \\
\end{array}
\right.
\end{equation}

Or equivalently, if we denote $z_\alpha=x_\alpha+i y_\alpha$ and $
\left\{
\begin{array}{c}
  \psi_1^\alpha=x_\alpha\circ\psi \\
   \psi_2^\alpha=y_\alpha\circ\psi\\
\end{array}
\right., $
 $\forall\alpha=1,...,p$,
we have:
\begin{equation}\label{}
\left\{
\begin{array}{ll}
\langle\nabla(x_\alpha\circ(\psi\circ\varphi)),
\nabla(x_\beta\circ(\psi\circ\varphi))\rangle -\\
\langle\nabla(y_\alpha\circ(\psi\circ\varphi)),
 \nabla(y_\beta\circ(\psi\circ\varphi))\rangle=0& \mbox{ a.e. in }X \\
\langle\nabla(y_\alpha\circ(\psi\circ\varphi)),
\nabla(x_\beta\circ(\psi\circ\varphi))\rangle+ \\
\langle\nabla(x_\alpha\circ(\psi\circ\varphi)),
\nabla(y_\beta\circ(\psi\circ\varphi))\rangle =0&\mbox{ a.e. in }X \\
\end{array}
\right.
\end{equation}
for every $\alpha,\beta =1,...,p$. Which means by definition, that
$\psi\circ\varphi$ is a PHWC map.

Conversely, suppose now that $\psi\circ\varphi$ is
pseudo-horizontally weakly conformal for any local holomorphic map
$\psi:N\rightarrow P$. Consider a local complex chart
$(z_\alpha)_{\alpha=1,...,p}$ in $P$. So the map
$\psi^\alpha\circ\varphi$, $\forall \alpha=1,...,p$ is
pseudo-horizontally weakly conformal (in the sense that we apply
(\ref{01}) for $v=w=z_\alpha$).

In order to prove that $\varphi$ is a PHWC map we shall use
Proposition \ref{prop:01'}. Let $u:N\rightarrow {\bf C}$ denote a
local holomorphic function on $N$. We associate to $u$ a new map:
$$
\begin{array}{cccc}
 \phi_u:& N &\rightarrow   & {\bf C}^p \\
   &x& \mapsto & (0,...,0,\underbrace{u(x)},0,...,0) \\
   & & &\alpha\\
\end{array}
$$

For any local holomorphic chart $\eta:P\rightarrow{\bf C}^p$, we
obtain a map $\psi_u:N\rightarrow P$, $\psi_u=\eta^{-1}\circ
\phi_u$ such that the $\alpha$'s coordinate of $\psi_u$ is exactly
the complex function $u:N\rightarrow {\bf C}$.

So we have proved that every local holomorphic function $u$ on $N$
can be obtained as a coordinate function of some local holomorphic
map $\psi_u:N\rightarrow P$.

We have supposed before, that for every $\psi:N\rightarrow P$ and
any local complex chart $(z_\alpha)_{\alpha=1,...,p}$ in $P$,
$$
\left\{
\begin{array}{cc}
\langle\nabla(\psi_1^\alpha\circ\varphi),\nabla(\psi_1^\alpha\circ\varphi)\rangle-
\langle\nabla(\psi_2^\alpha\circ\varphi),\nabla(\psi_2^\alpha\circ\varphi)\rangle=0&
\mbox {  a.e. in X}\\
\langle\nabla(\psi_2^\alpha\circ\varphi),\nabla(\psi_1^\alpha\circ\varphi)\rangle+
\langle\nabla(\psi_1^\alpha\circ\varphi),\nabla(\psi_2^\alpha\circ\varphi)\rangle=0&
\mbox {  a.e. in X} \\
\end{array}
\right.
$$

In particular, for $\psi=\psi_u$ ( $\psi^\alpha=u$), we obtain
$u\circ \varphi$ is pseudo-horizontally weakly conformal, for any
$u\in Holom (N)$. This implies (cf. Proposition \ref{prop:01'})
that $\varphi$ is a pseudo-horizontally weakly conformal map.
\endproof

\section{Pseudo harmonic morphism on Riemannian polyhedra.}

Similarly to the smooth case, if the target manifold is endowed
with a K\"ahler  structure one can enlarge the class of harmonic
morphisms on Riemannian polyhedra to the class of {\it pseudo
harmonic morphisms}.

Let $(X,g)$ denote an admissible Riemannian polyhedron and
$(N,J^N,g_N)$ a K\"ahler manifold without boundary.

\begin{definition}A map $\varphi:X\rightarrow N$ is called
{\rm pseudo harmonic morphism }(shortening PHM) if and only if
$\varphi$ is a harmonic map ( in the sense of Korevaar-Schoen {\rm
\cite {KS}} and Eells-Fuglede {\rm \cite {EF}}) and
pseudo-horizontally weakly conformal.
\end{definition}

Now, we will give a characterization of a pseudo harmonic morphism
in terms of the germs of the holomorphic functions on the target
K\"{a}hler manifold, and the harmonic structure (in sense of
Brelot, see ch 2 {\rm \cite {EF}}) on the domain admissible
Riemannian polyhedron.

\begin{thm}\label{thm:01}
 A continuous map
$\varphi:X\rightarrow N$ of class $W_{loc}^{1,2}(X,N)$ is a pseudo
harmonic morphism if and only if $\varphi$ pulls back local
complex-valued holomorphic functions on $N$ to harmonic functions
on $X$ (i.e. for any holomorphic function $ \psi:V\rightarrow {\bf
C}$ defined on a open subset $V$ of $N$ with $\varphi^{-1}(V)$
non-empty,  the composition
$\psi\circ\varphi:\varphi^{-1}(V)\rightarrow {\bf R}^2$ is
harmonic if and only if $\varphi$ is PHM).
\end{thm}

In order to prove the theorem, we shall need the following
elementary lemma:

\begin{lemma}\label{l-1}
Let $\varphi\in W_{loc}^{1,2}(X,N)$ and $\psi\in Holom (N)$. If
$\psi=\psi^1+i\psi^2$, then $\psi\circ\varphi\in
W_{loc}^{1,2}(X,{\bf R}^2)$ and moreover $\psi^j\circ\varphi\in
W_{loc}^{1,2}(X)$, $\forall j=1,2$.
\end{lemma}

\proof Let $\psi:N\rightarrow {\bf C}$ be a holomorphic function
and $K\subset N$ denote a compact subset in $N$. We say that
$\psi$ is uniformly Lipschitz in $K$ if there is a scalar
$\lambda$ (depending on $K$), such that $||\psi(p)-\psi(q)||\leq
\lambda d_N(p,q)$, $\forall p,q\in K$, where $||.||$ denotes the
usually norm in ${\bf C}$ and $d_N$ is the associated distance to
the  Riemannian metric $g_N$ on $N$. So $\psi$ is uniformly
Lipschitz on $K$.

We take $U$ a quasiopen set on $X$ which is relatively compact.
$\psi_{|\varphi(\overline{U})}$ is uniformly Lipschitz.

Now, following a result of Eells-Fuglede (see Corollary 9.1,
p.158, \cite{EF})
$$
E_U(\psi\circ\varphi)\leq
\lambda^2_{\varphi(\overline{U})}E_U(\varphi).
$$
It suffices to cover $X$ by a countable quasiopen sets $U$, such
that $E_U(\varphi)<\infty$ (such thing is possible because
$\varphi\in W_{loc}^{1,2}(X,N)$). We conclude that
$\psi\circ\varphi\in W_{loc}^{1,2}(X,{\bf R}^2)$.

For $\psi^j\circ\varphi$, $\forall j=1,2$, just remark that
$\psi^j$, $\forall j=1,2$ are locally Lipschitz, so by the same
argument used for proving that $\psi\circ\varphi\in
W_{loc}^{1,2}(X,{\bf R}^2)$, we have $\psi^j\circ\varphi\in
W_{loc}^{1,2}(X)$, $\forall j=1,2$.
\endproof

\begin{rmk}\label{rmk:th}
Let $f:N\rightarrow {\bf C}$ be a local holomorphic function
defined on a complex manifold $N$. Take $(z_1,...,z_n)$ a local
holomorphic coordinates in $N$ and denote $z_j=x_j+y_j, \forall
j=1,...,n$. If the manifold $N$ is K\"ahler and if we write
$f=f^1+if^2$, than we have the following equalities
("symmetries"):
$$
\frac{\partial^2 f^j}{\partial x_A\partial y_B}= \frac{\partial^2
f^j}{\partial x_B\partial y_A}, \forall j=1,2; \forall
A,B=1,...,n.
$$
\end{rmk}

\proofth{\bf 4.2:}

$"\Rightarrow"$:

 Let $\varphi:X\rightarrow N$ be a continuous map
of class $W_{loc}^{1,2}(X,N)$. Suppose that $\varphi$ is
pseudo-harmonic morphism and let $\psi:N\rightarrow {\bf C}$ be
any local holomorphic function, with $\psi=\psi^1+ i \psi^2$.
Remark that $\psi\circ\varphi$ is a continuous map ( as a
composite of two continuous maps). Moreover, by Lemma \ref{l-1},
$\psi\circ\varphi\in W_{loc}^{1,2}(X,{\bf R^2})$.

In order to prove that $\psi\circ\varphi$ is harmonic, (cf.
subsection 2.3.1) it suffices  to show that $\psi\circ\varphi$ is
weakly harmonic (as a map with values in ${\bf R}^2$).

 The fact that $\psi\circ\varphi$ is weakly
harmonic reads in the unique (conformal) chart $\eta:{\bf
R}^2\rightarrow {\bf R}^2$ and for any quasiopen set
$U\subset\varphi^{-1}({\bf R}^2)$ of compact closure in $X$,
\begin{equation}\label{13}
\begin{array}{ccc}
 &\int\limits_ U\langle\nabla f,\nabla(\psi^i\circ\varphi)\rangle d\mu_g=&\\
&=\int\limits_ U f.(^{{\bf
R}^2}\Gamma^i_{\alpha\beta}\circ\varphi)
\langle\nabla(\psi^\alpha\circ\varphi),\nabla(\psi^\beta\circ\varphi)\rangle d\mu_g,&\\
\end{array}
\end{equation}
 for $i=1,2$ and every bounded function $f\in
W_c^{1,2}(U)$.

In the case of ${\bf R}^2$ we know that $^{{\bf
R}^2}\Gamma^i_{\alpha\beta}\equiv0$, for any $\alpha,\beta=1,2$,
so the equations (\ref{13}) become:
\begin{equation}\label{14}
\int\limits_ U\langle\nabla f,\nabla(\psi^i\circ\varphi)\rangle
d\mu_g=0, \forall i=1,2.
\end{equation}
Let us  prove now (\ref{14}).

Take a real chart $(x_1,...,x_n,y_1,...,y_n)$ in $N$, such that
the complex associated chart $(z_j=x_j+ i y_j)_{j=1,n}$ is
holomorphic. Then for any chart domain $V\subset N$ and quasiopen
set $U\subset\varphi^{-1}(V)$ of compact closure in $X$, and for
any functions $v\in \mathcal{C}^2(V)$ and $f\in W_c^{1,2}(U)\cap
L^\infty(U)$ we have, (cf. \cite{EF} Remark 9.7):
\begin{equation}\label{15}
\int\limits_ U\langle\nabla f,\nabla(v\circ\varphi)\rangle d\mu_g=
\int\limits_ U\langle[(\partial_\alpha v)\circ\varphi]\nabla f,
\nabla\varphi^\alpha\rangle d\mu_g,
\end{equation}
where $\partial_\alpha$ denotes the $\alpha$'s partial derivative.

By partial integration,
\begin{equation}\label{16}
\begin{array}{cc}
\int\limits_ U\langle\nabla f,\nabla(v\circ\varphi)\rangle d\mu_g=
& \int\limits_ U\langle\nabla(f.[(\partial_\alpha
v)\circ\varphi]),
\nabla\varphi^\alpha\rangle d\mu_g -\\
  & \int\limits_ U f.[(\partial_\alpha\partial_\beta v)\circ\varphi]
  \langle\nabla \varphi^\alpha,
\nabla\varphi^\beta\rangle d\mu_g \\
\end{array}
\end{equation}
Recall that in the local coordinates $(x_1,...,x_n,y_1,...,y_n)$
in $N$, we have:
\begin{equation}\label{17}
v_{\alpha\beta}=\partial_\alpha\partial_\beta v-
^N\Gamma^k_{\alpha\beta}\partial_k v, \forall \alpha,\beta=1,2,
\end{equation}
(where $v_{\alpha\beta}$ are the second order covariant
derivatives of $v$). Inserting (\ref{17}) in the last integral on
the righthand side of the equation (\ref{16}) we obtain:
\begin{equation}\label{18}
\begin{array}{ccc}
 \int\limits_ U\langle\nabla f,\nabla(v\circ\varphi)\rangle d\mu_g&= &
  \int\limits_ U\langle\nabla(f.[(\partial_\alpha v)\circ\varphi]),
\nabla\varphi^\alpha\rangle d\mu_g -\\
  & &\int\limits_ U f.(v_{\alpha\beta}\circ\varphi)
\langle\nabla\varphi^\alpha,\nabla\varphi^\beta\rangle d\mu_g -\\
  &  &\int\limits_ U f.[(\partial_k v)\circ\varphi]
  (^N\Gamma^k_{\alpha\beta}\circ\varphi)
  \langle\nabla\varphi^\alpha,\nabla\varphi^\beta\rangle d\mu_g.\\
\end{array}
\end{equation}
$\varphi$ is supposed harmonic so it is weakly harmonic.
Consequently, the first and the third integral in the righthand
side of (\ref{18}) are equal, so:
$$
\int\limits_ U\langle\nabla f,\nabla(v\circ\varphi)\rangle d\mu_g=
-\int\limits_ U f.(v_{\alpha\beta}\circ\varphi)
\langle\nabla\varphi^\alpha,\nabla\varphi^\beta\rangle d\mu_g.
$$
Take now $v=\psi^i$. Then for $i=1,2$ we have:
$$
\int\limits_ U\langle\nabla f,\nabla(\psi^i\circ\varphi)\rangle
d\mu_g= -\int\limits_ U f.(\psi^i_{\alpha\beta}\circ\varphi)
\langle\nabla\varphi^\alpha,\nabla\varphi^\beta\rangle d\mu_g.
$$
We compute the righthand side of the last equality:
$$
\begin{array}{ccc}
 &\int\limits_ U f.(\psi^i_{\alpha\beta}\circ\varphi)
\langle\nabla\varphi^\alpha,\nabla\varphi^\beta\rangle d\mu_g  = & \\
 &=\int\limits_ U f.\sum\limits_{A,B=1}^n
(\frac{\partial^2\psi^i}{\partial x_A\partial x_B}\circ\varphi)
\langle\nabla(x_A\circ \varphi),\nabla(x_B\circ \varphi)\rangle d\mu_g+& \\
 &\int\limits_ U f.\sum\limits_{A,B=1}^n
(\frac{\partial^2\psi^i}{\partial x_A\partial y_B}\circ\varphi)
\langle\nabla(x_A\circ \varphi),\nabla(y_B\circ \varphi)\rangle d\mu_g+ &\\
 &\int\limits_ U f.\sum\limits_{A,B=1}^n
(\frac{\partial^2\psi^i}{\partial y_A\partial x_B}\circ\varphi)
\langle\nabla(y_A\circ \varphi),\nabla(x_B\circ \varphi)\rangle d\mu_g+& \\
&\int\limits_ U f. \sum\limits_{A,B=1}^n
(\frac{\partial^2\psi^i}{\partial y_A\partial y_B}\circ\varphi)
\langle\nabla(y_A\circ \varphi),\nabla(y_B\circ \varphi)\rangle d\mu_g. &\\
\end{array}
$$
For example, for $i=1$ we have:
$$
\begin{array}{ccc}
 &\int\limits_ U f.(\psi^1_{\alpha\beta}\circ\varphi)
\langle\nabla\varphi^\alpha,\nabla\varphi^\beta\rangle d\mu_g  = & \\
 & \int\limits_ U f.\sum\limits_{A,B=1}^n
(\frac{\partial^2\psi^1}{\partial x_A\partial x_B}\circ\varphi)
\langle\nabla(x_A\circ \varphi),\nabla(x_B\circ \varphi)\rangle d\mu_g+& \\
 &\int\limits_ U f. \sum\limits_{A,B=1}^n
(\frac{\partial^2\psi^1}{\partial x_A\partial y_B}\circ\varphi)
\langle\nabla(x_A\circ \varphi),\nabla(y_B\circ \varphi)\rangle d\mu_g+ &\\
 &\int\limits_ U f.\sum\limits_{A,B=1}^n
(\frac{\partial^2\psi^1}{\partial y_A\partial x_B}\circ\varphi)
\langle\nabla(y_A\circ \varphi),\nabla(x_B\circ \varphi)\rangle d\mu_g+& \\
&\int\limits_ U f.\sum\limits_{A,B=1}^n
(\frac{\partial^2\psi^1}{\partial y_A\partial y_B}\circ\varphi)
\langle\nabla(y_A\circ \varphi),\nabla(y_B\circ \varphi)\rangle d\mu_g. &\\
\end{array}
$$
Inserting  Cauchy-Riemann equations associated to $\psi_1$ and
$\psi_2$ (because $\psi$ is holomorphic) in the second and third
sums of the righthand side of the above equality, we obtain:
$$
\begin{array}{ccc}
 &\int\limits_ U f.(\psi^1_{\alpha\beta}\circ\varphi)
\langle\nabla\varphi^\alpha,\nabla\varphi^\beta\rangle d\mu_g  = & \\
 & \int\limits_ U f. \sum\limits_{A,B=1}^n
(\frac{\partial^2\psi^1}{\partial x_A\partial x_B}\circ\varphi)
\langle\nabla(x_A\circ \varphi),\nabla(x_B\circ \varphi)\rangle d\mu_g+& \\
 &\int\limits_ U f. \sum\limits_{A,B=1}^n
(-\frac{\partial^2\psi^2}{\partial x_A\partial x_B}\circ\varphi)
\langle\nabla(x_A\circ \varphi),\nabla(y_B\circ \varphi)\rangle d\mu_g+ &\\
 &\int\limits_ U f. \sum\limits_{A,B=1}^n
(\frac{\partial^2\psi^2}{\partial y_A\partial y_B}\circ\varphi)
\langle\nabla(y_A\circ \varphi),\nabla(x_B\circ \varphi)\rangle d\mu_g+& \\
&\int\limits_ U f. \sum\limits_{A,B=1}^n
(\frac{\partial^2\psi^1}{\partial y_A\partial y_B}\circ\varphi)
\langle\nabla(y_A\circ \varphi),\nabla(y_B\circ \varphi)\rangle d\mu_g &\\
\end{array}$$

 By easy computations we obtain:
$$
\begin{array}{ccc}
&\int\limits_ U f.(\psi^1_{\alpha\beta}\circ\varphi)
\langle\nabla\varphi^\alpha,\nabla\varphi^\beta\rangle d\mu_g=
\int\limits_ U f. \sum\limits_{A,B=1}^n
[(\frac{\partial^2\psi^1}{\partial x_A\partial x_B}\circ\varphi)-
(\frac{\partial^2\psi^1}{\partial y_A\partial y_B}\circ\varphi)]&\\
&[\langle\nabla(x_A\circ \varphi),\nabla(x_B\circ \varphi)\rangle-
\langle\nabla(y_A\circ \varphi),\nabla(y_B\circ \varphi)\rangle]d\mu_g-&\\
&\int\limits_ U f. \sum\limits_{A,B=1}^n
[(\frac{\partial^2\psi^2}{\partial x_A\partial x_B}\circ\varphi)-
(\frac{\partial^2\psi^2}{\partial y_A\partial y_B}\circ\varphi)]&\\
&[\langle\nabla(x_A\circ \varphi),\nabla(y_B\circ \varphi)\rangle+
\langle\nabla(y_A\circ \varphi),\nabla(x_B\circ \varphi)\rangle]d\mu_g+&\\
&\int\limits_ U f. \sum\limits_{A,B=1}^n
[(\frac{\partial^2\psi^1}{\partial x_A\partial x_B}\circ\varphi)
\langle\nabla(y_A\circ \varphi),\nabla(y_B\circ \varphi)\rangle+&\\
&(\frac{\partial^2\psi^1}{\partial y_A\partial y_B}\circ\varphi)
\langle\nabla(x_A\circ \varphi),\nabla(x_B\circ \varphi)\rangle]d\mu_g+&\\
&\int\limits_ U f. \sum\limits_{A,B=1}^n
[(\frac{\partial^2\psi^2}{\partial x_A\partial x_B}\circ\varphi)
\langle\nabla(y_A\circ \varphi),\nabla(x_B\circ \varphi)\rangle-&\\
&(\frac{\partial^2\psi^2}{\partial y_A\partial y_B}\circ\varphi)
\langle\nabla(x_A\circ \varphi),\nabla(y_B\circ \varphi)\rangle]d\mu_g.&\\
\end{array}
$$
Because the map $\varphi$ is supposed pseudo-horizontally weakly
conformal,  the first and the second integral of the righthand
side in the last equality are zero, so:
$$
\begin{array}{ccc}
&\int\limits_ U f.(\psi^1_{\alpha\beta}\circ\varphi)
\langle\nabla\varphi^\alpha,\nabla\varphi^\beta\rangle d\mu_g = \\
 &\int\limits_ U f. \sum\limits_{A,B=1}^n
 [(\frac{\partial^2\psi^1}{\partial x_A\partial x_B}\circ\varphi)+
 (\frac{\partial^2\psi^1}{\partial y_A\partial y_B}\circ\varphi)]
  \langle\nabla(x_A\circ \varphi),\nabla(x_B\circ \varphi)\rangle d\mu_g+&  \\
 & +\int\limits_ U f. \sum\limits_{A,B=1}^n
 [(\frac{\partial^2\psi^2}{\partial x_A\partial x_B}\circ\varphi)+
 (\frac{\partial^2\psi^2}{\partial y_A\partial y_B}\circ\varphi)]
 \langle\nabla(x_A\circ \varphi),\nabla(y_B\circ \varphi)\rangle d\mu_g & \\
\end{array}
$$
Now, by Remark (\ref{rmk:th}) the last two sums are zero. So we
obtain:
$$
\int\limits_ U\langle\nabla f,\nabla(\psi^1\circ\varphi)\rangle
d\mu_g=0.
$$
By a similar computation we can also prove:
$$
\int\limits_ U\langle\nabla f,\nabla(\psi^2\circ\varphi)\rangle
d\mu_g=0.
$$
Thus $\psi\circ\varphi$ is weakly harmonic and so, harmonic (cf.
subsection 2.3.1).

$ "\Leftarrow"$:

Conversely, suppose $\psi\circ\varphi:X\rightarrow {\bf R}^2$ is
harmonic  for any $\psi\in Holom(N)$. It is known that the
harmonicity of $\psi\circ\varphi$ is equivalent (cf. subsection
2.3.1 and Lemma \ref{l-1}) to the weak harmonicity of
$\psi\circ\varphi$. So, for a given chart $\eta:{\bf
R}^2\rightarrow {\bf R}^2$, and for any quasiopen set $U\subset
\varphi^{-1}({\bf R}^2)$ of compact closure in $X$, the weak
harmonicity  reads:
\begin{equation}\label{19}
\int\limits_U\langle\nabla f,\nabla(\psi^i\circ\varphi)\rangle
d\mu_g=0,
\end{equation}
 for $i=1,2$ and every function $f\in W_c^{1,2}(U)$,
 where $\psi=\psi_1 +i \psi_2$:.

Now, choose local complex holomorphic coordinates in $N$,
$(z_1,...,z_n)$, where $z_k=x_k+iy_k$, $\forall k=1,...,n$. Take
$\psi=z_k$, $\forall k=1,...,n$, and denote
$$
\begin{array}{cccccc}
  x_1,     & ..., &x_n,      &  y_1,      & ..., & y_n \\
\parallel  &      &\parallel &  \parallel &      & \parallel \\
\xi_1      &      &\xi_n     &  \xi_{n+1} &      & \xi_{2n}\\
\end{array}
$$
by a generic term $\xi^k$, $\forall k=1,...,2n$. In particular,
the equation (\ref{19}), for the coordinates functions $\xi^k$,
reads:
\begin{equation}\label{20}
\int\limits_U\langle\nabla f,\nabla(\xi^k\circ\varphi)\rangle
d\mu_g=0, \forall k=1,...,2n.
\end{equation}
For the domain chart $V$ of $(z_k)$, any quasiopen set
$U\subset\varphi^{-1}(V)$ of compact closure in $X$ and for $f\in
W^{1,2}_c(U)$ we have:
$$
\begin{array}{cc}
  0 &= \int\limits_U\langle\nabla f,\nabla(\xi^k\circ\varphi)\rangle d\mu_g=
\int\limits_U\langle\nabla
(f.[(\partial_\alpha\xi^k)\circ\varphi]),
\nabla\varphi^\alpha\rangle d\mu_g-\\
   &\int\limits_U f.(\xi^k_{\alpha\beta}\circ\varphi)
   \langle\nabla\varphi^\alpha,\nabla\varphi^\beta\rangle d\mu_g-\\
   &\int\limits_U f.[(\partial_\gamma\xi^k)\circ\varphi]
   (^N\Gamma ^\gamma_{\alpha\beta}\circ\varphi)
   \langle\nabla\varphi^\alpha,\nabla\varphi^\beta\rangle d\mu_g,\\
   &\forall \alpha,\beta,\gamma=1,...,2n.\\
\end{array}
$$
Taking into account that: $\forall k=1,...,2n$ and $\forall
\alpha,\beta=1,...,2n$,
$$\xi^k_{\alpha\beta}=0 \mbox{ and }
\partial_\alpha\xi^k=
\left\{
\begin{array}{cc}
  0,& \mbox {if } \alpha=k \\
  1, & \mbox {if } \alpha\neq k\\
\end{array}
\right.,
$$
we obtain:
$$
\begin{array}{cc}
  &0 = \int\limits_U\langle\nabla f,\nabla(\xi^k\circ\varphi)\rangle d\mu_g=\\
  &= \int\limits_U\langle\nabla f,\nabla \varphi^k\rangle d\mu_g-
  \int\limits_U f.
   (^N\Gamma ^k_{\alpha\beta}\circ\varphi)
   \langle\nabla\varphi^\alpha,\nabla\varphi^\beta\rangle d\mu_g,\\
\end{array}
$$

So, $ \forall \alpha,\beta,k=1,...,2n$,
$$
\int\limits_U\langle\nabla f,\nabla\varphi^k\rangle d\mu_g=
\int\limits_U f.
   (^N\Gamma ^k_{\alpha\beta}\circ\varphi)
   \langle\nabla\varphi^\alpha,\nabla\varphi^\beta\rangle d\mu_g.
 $$
 This means that $\varphi$ is weakly harmonic.
 But $\varphi$ is a continuous map of class
 $W_{loc}^{1,2}(X,N)$, so it is a harmonic map (cf. subsection 2.3.1).

 Now, for any $v:N\rightarrow {\bf R}$
 such that $v\in\mathcal{C}^2(V)$, where $V$
 is a domain chart on $N$, we have:
 \begin{equation}\label{21}
 \begin{array}{cc}
 \int\limits_U\langle\nabla f,\nabla(v\circ\varphi)\rangle d\mu_g =&
 \int\limits_U\langle\nabla (f[(\partial_kv)\circ\varphi]),
 \nabla\varphi^k\rangle d\mu_g - \\
  & \int\limits_U f(v_{\alpha\beta}\circ\varphi)
  \langle\nabla\varphi^\alpha,\nabla\varphi^\beta\rangle d\mu_g-\\
  &\int\limits_U f[(\partial_kv)\circ\varphi](^N\Gamma^k_{\alpha\beta}\circ\varphi)
  \langle\nabla\varphi^\alpha,\nabla\varphi^\beta\rangle d\mu_g,\\
 \end{array}
\end{equation}
for $\alpha,\beta=1,...,2n$, where $U\subset\varphi^{-1}(V)$ is
quasiopen with compact closure in $X$.

In particular, for any holomorphic chart:
$$
\begin{array}{cccc}
\eta:V\subset N & \rightarrow  & {\bf C}^n & \\
             p  & \mapsto       & (z_1,...,z_n), & z_j=x_j+iy_j, \\
\end{array}
$$
and for any local holomorphic function $\psi:N\rightarrow {\bf
C}$, $\psi=\psi^1+i\psi^2$, denoting $x_1,...,x_n,y_1,...,y_n$ by
$\xi^\gamma$ , $\forall \gamma=1,...,2n$, we apply equation (\ref
{21}) to $\psi^1$ and $\psi^2$ respectively:
$$
\begin{array}{ccc}
 0=&
 \int\limits_U\langle\nabla (f[(\partial_k\psi^i)\circ\varphi]),
 \nabla(\xi^k\circ\varphi)\rangle d\mu_g - &\\
  & \int\limits_U f(\psi^i_{\alpha\beta}\circ\varphi)
  \langle\nabla(\xi^\alpha\circ\varphi),\nabla(\xi^\beta\circ\varphi)\rangle d\mu_g -&\\
  &\int\limits_U f[(\partial_k\psi^i)\circ\varphi](^N\Gamma^k_{\alpha\beta}\circ\varphi)
  \langle\nabla(\xi^\alpha\circ\varphi),\nabla(\xi^\beta\circ\varphi)\rangle d\mu_g,
  &\forall i=1,2.\\
\end{array}
$$
In the above equality, the first and the last integral are equal
(because $\varphi$ is harmonic), so we have:
$$
0=\int\limits_U f(\psi^i_{\alpha\beta}\circ\varphi)
  \langle\nabla(\xi^\alpha\circ\varphi),\nabla(\xi^\beta\circ\varphi)\rangle d\mu_g.
$$

In the proof of the "if" part we have obtained:
\begin{equation}\label{22}
\begin{array}{ccc}
&0=\int\limits_U f(\psi^1_{\alpha\beta}\circ\varphi)
  \langle\nabla(\xi^\alpha\circ\varphi),\nabla(\xi^\beta\circ\varphi)\rangle d\mu_g=&\\
&=\int\limits_ U f.\sum\limits_{A,B=1}^n
[(\frac{\partial^2\psi^1}{\partial x_A\partial x_B}\circ\varphi)-
(\frac{\partial^2\psi^1}{\partial y_A\partial y_B}\circ\varphi)]&\\
&[\langle\nabla(x_A\circ \varphi),\nabla(x_B\circ \varphi)\rangle-
\langle\nabla(y_A\circ \varphi),\nabla(y_B\circ \varphi)\rangle]d\mu_g -&\\
&\int\limits_ U f.\sum\limits_{A,B=1}^n
[(\frac{\partial^2\psi^2}{\partial x_A\partial x_B}\circ\varphi)-
(\frac{\partial^2\psi^2}{\partial y_A\partial y_B}\circ\varphi)]&\\
&[\langle\nabla(x_A\circ \varphi),\nabla(y_B\circ \varphi)\rangle+
\langle\nabla(y_A\circ \varphi),\nabla(x_B\circ \varphi)\rangle]d\mu_g+&\\
&\int\limits_ U f. \sum\limits_{A,B=1}^n
[(\frac{\partial^2\psi^1}{\partial x_A\partial x_B}\circ\varphi)
\langle\nabla(y_A\circ \varphi),\nabla(y_B\circ \varphi)\rangle+&\\
&(\frac{\partial^2\psi^1}{\partial y_A\partial y_B}\circ\varphi)
\langle\nabla(x_A\circ \varphi),\nabla(x_B\circ \varphi)\rangle]d\mu_g+&\\
&\int\limits_ U f. \sum\limits_{A,B=1}^n
[(\frac{\partial^2\psi^2}{\partial x_A\partial x_B}\circ\varphi)
\langle\nabla(y_A\circ \varphi),\nabla(x_B\circ \varphi)\rangle -&\\
&(\frac{\partial^2\psi^2}{\partial y_A\partial y_B}\circ\varphi)
\langle\nabla(x_A\circ \varphi),\nabla(y_B\circ \varphi)\rangle]d\mu_g.&\\
\end{array}
\end{equation}
Choosing particular holomorphic functions as: $\psi=z_Az_B$ and
$\psi=iz_Az_B$, $\forall A,B=1,...,n$, in equation (\ref{22}), we
obtain the pseudo-horizontally weakly conformal conditions (see
Definition \ref{def:01}).
\endproof

A straightforward result deduced from the Theorem \ref{thm:01} is
the folowing:

\begin{cor}\label{cor:01}
Let $\varphi:X\rightarrow N$ be a pseudo-horizontally weakly
conformal map from a Riemannian polyhedron into a K\"ahler
manifold. Then $\varphi$ is harmonic if and only if the components
$\varphi^k=z_k\circ\varphi$ in terms of any holomorphic
coordinates $(z_k)_{k=1,...,n}$ in $N$, are harmonic maps
($\varphi^k$ is understood as a map into ${\bf R}^2$).
\end{cor}

\proof The "only if part" is already proved through the proof of
the Theorem \ref{thm:01}.

The "if part": Suppose that the map $\varphi$ is harmonic. So it
is pseudo harmonic morphism. By Theorem \ref{thm:01}, $\varphi$
pulls back local holomorphic functions to local harmonic maps
which applies to any holomorphic coordinates in $N$.
\endproof

We turn our attention  to the relation between  pseudo harmonic
morphism on Riemannian polyhedron and local holomorphic maps
between K\"{a}hler manifolds. More precisely, we have:

\begin{thm}
Let $(X,g)$ be a Riemannian polyhedron and $(N,J^N,g_N)$,
$(P,J^P,g_P)$ be two K\"ahler manifolds with $dim_{\bf C}N=n$ and
$dim_{\bf C}P=p$. A continuous map $\varphi:X\rightarrow N$, of
class $W_{loc}^{1,2}(X,N)$, is pseudo harmonic morphism if and
only if $\psi\circ\varphi:X\rightarrow P$ is a (local) pseudo
harmonic morphism for all local holomorphic maps
$\psi:N\rightarrow P$.
\end{thm}

\proof $"\Rightarrow":$

Suppose that $\varphi:X\rightarrow N$ is a pseudo harmonic
morphism and let $\psi:X\rightarrow P$ be any local holomorphic
map.

By Proposition \ref{prop:01}, $\psi\circ\varphi$ is a
pseudo-horizontally weakly conformal map.

 Let $(z_1,...,z_p)$ be local holomorphic coordinates in $P$, and set
$\psi^j=z_j\circ\psi$. On one hand, the map $\varphi$ is supposed
PHM and obviously the complex functions $\psi^j:P\rightarrow {\bf
C}$, $\forall j=1,...,p$, are holomorphic so,  by Theorem
\ref{thm:01}, $\psi^j\circ\varphi$ are harmonic, $\forall
j=1,...,p$.

On the other hand, $\psi\circ\varphi$ is pseudo-horizontally
weakly conformal. Thus by Corollary \ref{cor:01}, the map
$\psi\circ\varphi$ is (local) harmonic.

$"\Leftarrow":$

Suppose now that $\psi\circ\varphi:X\rightarrow P$ is a (local)
pseudo harmonic morphism, for any local holomorphic map
$\psi:N\rightarrow P$. By Proposition \ref{prop:01}, the map
$\varphi$ is pseudo-horizontally weakly conformal.

As we have already done in the proof of Proposition \ref{prop:01},
 every local holomorphic function $v:N\rightarrow {\bf C}$ can be
 obtained  as a coordinate of some local holomorphic map $\psi_v:N\rightarrow P$.

$\psi_v\circ\varphi$ is local a pseudo harmonic morphism. By
Corollary \ref{cor:01}, $v\circ\varphi$ is harmonic (because it is
a coordinate of $\psi_v\circ\varphi$). Therefore, we have shown
that for any local holomorphic function $v:N\rightarrow {\bf C}$,
the map $v\circ\varphi$ is harmonic. So by Theorem \ref{thm:01} we
conclude that $\varphi$ is a pseudo harmonic morphism.
\endproof

 The following proposition will play an important roll in the next
 section for constructing examples of PHM on admissible Riemannian
 polyhedra.

\begin{prop}\label{prop:fin}
Let $(X,g)$ and $(Y,h)$ be two  admissible Riemannian polyhedra,
and $(N,J^N,g_N)$ a K\"ahler manifold without boundary of complex
dimension $n$. Let $\varphi:Y\rightarrow N$ be a continuous map of
class $W^{1,2}_{loc}(Y,N)$, $\pi:X\rightarrow Y$ a proper,
surjective, continuous map of class $W^{1,2}_{loc}(X,Y)$ and
$\widetilde{\varphi}=\varphi\circ \pi$. If $\pi$ is a harmonic
morphism, then $\varphi$ is pseudo harmonic morphism if and only
if $\widetilde{\varphi}$ is pseudo harmonic morphism.
\end{prop}

\proof Suppose that $\widetilde{\varphi} : X\rightarrow N$ is a
pseudo harmonic morphism; by Theorem \ref{thm:01}, it is
equivalent to the fact that $\widetilde{\varphi}$ pulls back local
holomorphic functions on $N$ to local harmonic functions on $X$.
So, for any $\psi\in Holom(N)$, $\psi=\psi_1+i \psi_2$, the map
$\psi\circ \widetilde{\varphi} \in W^{1,2}_{loc}(X,{\bf R}^2)$ is
locally harmonic. But this can also reads: for any $\psi\in
Holom(N)$, $\psi=\psi_1+i \psi_2$, the functions $\psi_1\circ
\widetilde{\varphi}$ and $\psi_2\circ \widetilde{\varphi}$ are
locally harmonic (because $\psi\circ \widetilde{\varphi}$ is
continuous and the Christoffel symbols relative to the fixed chart
of ${\bf R}^2$ are all zero [cf. subsection 2.3.1]); or, also for
any $\psi\in Holom(N)$, with $\psi=\psi_1+i \psi_2$, the functions
$(\psi_1\circ \varphi)\circ \pi$ and $(\psi_2\circ \varphi)\circ
\pi$ are locally harmonic.

On the other hand, the map $\pi$ is a harmonic morphism and it is
supposed surjective and proper so, by subsection 2.3.2 (cf.
\cite{EF}, Theorem 13.1), for every open set $V\subset Y$ and any
function $v:V\rightarrow {\bf R}$, the function $v$ is harmonic if
and only if $v\circ\pi$ is harmonic. In particular this fact
applies to the functions $(\psi_1\circ \varphi)$ and $(\psi_2\circ
\varphi)$; in other words, for any local holomorphic function
$\psi=\psi_1+i \psi_2$ on $N$, the functions $(\psi_1\circ
\varphi)\circ \pi$ respectively $(\psi_2\circ \varphi)\circ \pi$
are locally harmonic if and only if the functions $(\psi_1\circ
\varphi)$ respectively $(\psi_2\circ \varphi)$ are locally
harmonic.

The last assertion is equivalent to the fact that, for any
$\psi\in Holom(N)$, the map $\psi\circ\varphi$ is locally
harmonic. But by Theorem \ref{thm:01}, this means (iff) that
$\varphi:Y\rightarrow N$ is a pseudo harmonic morphism. This ends
the proof of the proposition.

\endproof

\section{Some examples.}

In this short section we will offer some examples of pseudo
harmonic morphisms on Riemannian polyhedra.

As we have shown in Proposition \ref{prop:00},  every horizontally
weakly conformal map from a Riemannian admissible polyhedron into
a Hermitian manifold, is pseudo horizontally weakly conformal. So
every harmonic morphism into a K\"ahler manifold is pseudo
harmonic morphism. For other nontrivial examples, when the source
polyhedra are smooth Riemannian manifolds, see for example
\cite{AAB}.

Thanks to Proposition \ref{prop:fin}, we can derive a several
non-obvious examples of pseudo harmonic morphisms on (singular)
Riemannian polyhedra.

It is known (see Example 8.12, \cite{EF}) that given a $K$ compact
group of isometries of a complete smooth Riemannian manifold $M$,
and $\pi:M\rightarrow M/K$ the projection onto the orbit space
$M/K$ (with the quotient topology), there is a smooth
triangulation of $M$ for which $\pi$ induces an admissible
Riemannian polyhedral structure on $M/K$. The associated intrinsic
distance $d_{M/K}(y_1,y_2)$ between elements $y_1$, $y_2$ of
$M/K$, equals the intrinsic distance in $M$ between the
corresponding compact orbits $\pi^{-1}(y_1)$ and $\pi^{-1}(y_2)$,
and the polyhedron $M/K$ with the distance $d_{M/K}$ is a geodesic
space. The polyhedral structure determines a Brelot harmonic sheaf
$\mathcal{H}_{M/K}$ on $M/K$ (cf. Theorem 7.1, \cite{EF}).
Moreover, $\pi$ caries the Brelot harmonic sheaf $\mathcal{H}_M$
onto a Brelot harmonic sheaf
$\mathcal{H}'_{M/K}=\pi_{*}\mathcal{H}_M$ on $M/K$ and
$\pi:(M,\mathcal{H}_M)\rightarrow(M/K,\mathcal{H}'_{M/K})$ becomes
a harmonic morphism, surjective and proper.

This construction can be applied to Riemannian orbifolds (cf.
Subexample 8.13(ii), \cite{EF}) as follows:

Let $M$ be a Riemannian manifold, and $\mathcal{S}^r$ the symmetry
group on $r$ factors.

Denote:
$$
\begin{array}{cc}
\mathcal{S}^r\mathcal{M}:=&(\underbrace{M\times M \times ...\times M})/\mathcal{S}^r\\
                          &r \mbox{ -times}                                         \\
\end{array}
$$
the $r$-fold symmetric power of the manifold $M$.

The compact group $\mathcal{S}^r$ acts isometrically, so
$\mathcal{S}^r\mathcal{M}$ becomes a Riemannian orbifold (singular
if the dimension of $M\geq3$). Thus, as above, we obtain a proper,
surjective, harmonic morphism from $M\times M \times ...\times M$
to $\mathcal{S}^r\mathcal{M}$.

In the particular case when $M={\bf C}^{k+s}$, with $2(k+s)\geq3$,
following \cite{AAB}, one  can construct a pseudo harmonic
morphism:
$$
\begin{array}{ccccc}
\eta&:{\bf C}^k\times{\bf C}^s&\rightarrow&{\bf C}^r& \mbox {  given by}\\
&(u,v)&\mapsto&\left(\frac{F_1(u)P_1(\overline{v})}{G_1(u)Q_1(\overline{v})},...,
\frac{F_r(u)P_r(\overline{v})}{G_r(u)Q_r(\overline{v})}\right),&\\
\end{array}
$$
where $F_1$,...,$F_r$,$G_1$,...,$G_r$ are homogenous polynomials
on ${\bf C}^k$ and  $P_1$,...,$P_r$,$Q_1$,...,$Q_r$ are homogenous
polynomials on ${\bf C}^s$, all having the same degree. We know
that the sum of two PHM is also a PHM, so we define:
$$
\widetilde{\varphi}:{\bf C}^{k+s}\times{\bf C}^{k+s}\rightarrow
{\bf C}^r, \mbox{   }
\widetilde{\varphi}(\widetilde{u},\widetilde{v}):=\eta(\widetilde{u})+\eta(\widetilde{v}),
$$
for any $\widetilde{u},\widetilde{v}\in {\bf C}^{k+s}$.

Using the harmonic morphism ${\bf C}^{k+s}\times{\bf
C}^{k+s}\rightarrow \mathcal{S}^2\mathcal {C}^{k+s}$, the map
$\widetilde{\varphi}$ factors through a PHM from
$\mathcal{S}^2\mathcal{C}^{k+s}$ to ${\bf C}^r$.

Moreover, using the harmonic morphism ${\bf C}^{k+s}\rightarrow
{\bf C}P^{k+s-1} $, the map $\eta$ factors (see \cite{AAB})
through a PHM from ${\bf C}P^{k+s-1}$ to ${\bf C}^r$ which is
neither holomorphic nor antiholomorphic. Now using this map and
apply the same arguments as before, we  get a pseudo harmonic
morphism from $\mathcal{S}^2\mathcal{CP}^{k+s-1} $ to  ${\bf
C}^r$.

\bigskip

\noindent {\em Acknowledgments.} We thank the ICTP Trieste for
hospitality during this work.

\end{document}